\documentclass[11pt,reqno]{amsart}

\usepackage{amssymb,amsmath,amsthm,graphics,epsfig}
\usepackage{amsfonts}
\usepackage{graphicx,cite}
\usepackage{amssymb}
\usepackage[pagewise]{lineno}
\usepackage{mathrsfs}
\usepackage[colorlinks=true, pdfstartview=FitV, linkcolor=blue, citecolor=blue, urlcolor=blue]{hyperref}
\usepackage{booktabs}
\usepackage{graphicx}
\usepackage{tikz}
\usepackage{subfigure}
\usepackage{subcaption}
\usepackage{amsthm}
\usepackage{multicol}
\usepackage{cases}
\usepackage{booktabs}
\usepackage{subcaption}

\usepackage{arydshln}
\usepackage[colorinlistoftodos]{todonotes}
\usepackage[normalem]{ulem}
\usepackage{color}
\usepackage{amsaddr} 
\usepackage[margin=2.7cm]{geometry}

\usepackage{lineno}

\newtheorem{Theorem}{Theorem}[section]

\usepackage{float}
\theoremstyle{definition}

\newtheorem{Definition}{Definition}

\definecolor{ColorEdward}{rgb}{0.4,0.6,0.7}

\definecolor{ColorFrancisco}{rgb}{0.956,0.878,0.5}

\title[Dynamics of the Three Fixed Centers Problem]{Stability, Periodic Orbits and KAM Tori in the Dynamics of the Three Fixed Centers Problem}

\author[E. Turner, F. Crespo, J. Vidarte, J. Villafañe and J. Zapata]{Edward A. Turner$^{1,2,*}$, Francisco Crespo$^{2}$, Jhon Vidarte$^{3}$\\ Jersson Villafañe$^2$ \MakeLowercase{and} Jorge Zapata$^4$}

\address{$^{1}${Universidad de Concepción, Departamento de Ciencias Básicas, Campus Los Ángeles,\\ {\small \itshape{Av. J.A. Coloma 0201, Los Ángeles, Chile}}}}
	
\address{$^{2}${Universidad del Bío-Bío, Grupo de investigaci\'on en Sistemas Din\'amicos\\ y Aplicaciones (GISDA), {\small \itshape{Casilla 5-C, Concepci\'on, Chile}}}}
\address{$^{3}$Universidad Católica de la Santísima Concepción, Departamento de Matemática\\ y Física,
{\small \itshape{Casilla 297, Concepción, Chile}}}
\address{$^{4}$Universidad San Sebastián, Facultad de Educación, Programa de Formaci\'on Pedag\'ogica\\ para Licenciados y/o Profesionales {\small \itshape{Lientur 1457, Concepción, Chile }}}

\thanks{$^{*}$Corresponding author: \texttt{edward.turner1901@alumnos.ubiobio.cl}}

	

\begin{document}
\maketitle	
\begin{abstract}
We investigate the motion in space of an infinitesimal particle in the gravitational field generated by three primary bodies positioned at the vertices of a fixed equilateral triangle. We assume that the distances between the primaries are small compared to their separation from the particle. By applying a Lie-Deprit normalization, we simplify the Hamiltonian, relegating both the mean anomaly and the argument of periapisis to third-order terms or higher.
After reducing out the symmetries associated with the Kepler flow and the central action of the angular momentum, we examine the relative equilibria in the first and second reduced spaces. 
We are able to identify the conditions for the existence of circular periodic orbits and KAM tori, thus providing insight into the system's long-term stability and dynamic structure.
\end{abstract}


\section{Introduction}
In 1760, Euler was the first to examine the $n$-center problem for the case $n=2$, which has since become one of the most renowned integrable problems in classical mechanics. This problem involves the motion of a test particle influenced by the gravitational field of two fixed centers. Jacobi later demonstrated that the equations governing this motion could be solved using elliptic functions \cite{Jacobi1969}. However, for $n \geq 3$, the problem becomes analytically non-integrable \cite{Bolotin1984(2), SHIBAYAMA20182461}. This system has served as a significant model in both macroscopic \cite{Chen2020, Holger2004,Coulson1967} and microscopic \cite{ Baber_Hassé_1935, Klein, Athavan2001(3)} contexts.  In molecular physics, the $n$-centers correspond to atomic nuclei, with the test particle representing an electron, providing a simple model for diatomic molecules. In celestial mechanics, it offers a simplified model of the motion of a fast-moving body influenced by $n$ fixed centers of attraction. This arises as an approximation of the more complex $n+1$-body problem, where one body moves much faster than the others, allowing for a reduction in complexity.

The dynamics of the particle are governed by the Hamiltonian function $\mathcal{H}: T^*\mathcal{N} \to \mathbb{R}$ on the phase space $T^*\mathcal{N}$, where the configuration space $\mathcal{N} = \mathbb{R}^3 - \{{\xi}_1, \dots, {\xi}_n\}$ excludes the fixed centers of attraction. The Hamiltonian is given by
\begin{equation}\label{eq:FullHam1}
\mathcal{H}({q}, {p}) = \frac{1}{2} \Vert {p} \Vert^2 - V\left({q}\right), \quad V({q}) = \sum_{i=1}^n \frac{\mu_i}{\lVert {q} - {\xi}_i \rVert},
\end{equation}
where $\mu_i$ represents the mass of the $i$-th particle, defined as $\mu_i = \mathcal{G} m_i \in \mathbb{R}$.

A large number of studies have focused on the case $n=2$. For example, Charlier \cite{Charlier1902} examined all possible motions in the planar two-center problem, which were later systematically analyzed by Deprit \cite{Deprit1962}. Waalkens et al. \cite{Holger2004} presented a global view of how the phase space is foliated by invariant tori, identifying topological changes or bifurcations involving special tori of lower dimension. Biscani and Izzo \cite{Biscani} provided a comprehensive, explicit, and complete solution to the three-dimensional problem of two fixed centers using Weierstrass elliptic and related functions.


Our work explores a specific system with $n=3$, where each center is located at the vertices of an equilateral triangle with equal masses. We study the dynamics of a particle that is far from the centers compared to the length of the edges, where the centers remain fixed. 

Later, we employ the Delaunay chart and introduce an artificial symmetry in the mean anomaly through a second-order Lie-Deprit normalization. In addition to eliminating the mean anomaly, this process relegates the appearance of the argument of the periapsis to higher-order terms, which enables a reduction process in two stages: the first reduction corresponds to eliminating the mean anomaly, while the second takes advantage of the absence of the argument of periapsis.

In the first reduced space \cite{Meyer1973, Marsden1974} and using Hamiltonian averaging theory \cite{Moser1970}, we identify six families of periodic solutions at each energy level and their stability. 
Additionally, we determine the existence of several families of 3-dimensional KAM tori in the first and second reduced spaces.

The paper is organized as follows: In Section~\ref{sec:2}, we formulate the problem and introduce the artificial symmetry associated with the Keplerian energy. In Section~\ref{sec:Reducido1}, we determine the existence of periodic orbits and KAM tori in the first reduced space. In Section~\ref{sec:Reducido2}, we identify new families of KAM tori from the twice-reduced space.
Finally, we include an appendix containing avereging and KAM theory.

\section{Dynamics of a particle far from a triple-star system}\label{sec:2}
In our formulation, we assume units such that the parameters $\mu_i$ are all equal to $1/3$, and the coordinates of the centers are given by 
\begin{equation}\label{eq:equilatero}
 \xi_1 = (0,0,0), \quad \xi_2 = (\varepsilon,0,0), \quad \xi_3 =(\varepsilon/2,\varepsilon\sqrt{3}/2,0),
 \end{equation} 
 where $\varepsilon > 0$ is a small parameter reflecting the assumption that the particle is moving far from the origin relative to the distances between the primaries and the origin. The phase space becomes $T^*\mathcal N$ with $\mathcal N=\mathbb R^3-\{\xi_1,\xi_2,\xi_3\}$ and the Hamiltonian function in Cartesian coordinates is defined by
\begin{equation}\label{eq:Ham3centros}
\mathcal H=\dfrac{1}{2}\Vert  p\Vert^2-\dfrac{1}{3}\left[ \dfrac{1}{\Vert q\Vert}+\dfrac{1}{\sqrt{(q_1-\varepsilon)^2+q_2^2+q_3^2}}+\dfrac{1}{\sqrt{(q_1-\varepsilon/2)^2+(q_2-\sqrt{3}\varepsilon/2)^2+q_3^2}}\right],
\end{equation}
where $q=(q_1,q_2,q_3)$ and $ p=(p_1,p_2,p_3).$
Then we consider the following expansion in powers of $\varepsilon$ up to the second order 
$
\mathcal H(q,p)={\mathcal H}_0+\varepsilon{\mathcal H}_1+(\varepsilon^2/2){\mathcal H}_2+O(\varepsilon^3),
$
where,
\begin{eqnarray}\label{eq:HamPertur}
 {\mathcal H}_0=\dfrac{1}{2}\Vert p\Vert^2-\dfrac{1}{\Vert q\Vert},
\quad {\mathcal H}_1=-\frac{3 q_1+\sqrt{3} q_2}{6 \left(q_1^2+q_2^2+q_3^2\right)^{3/2}},
\quad {\mathcal H}_{2}=-\frac{7 q_1^2+6 \sqrt{3} q_1 q_2+q_2^2-8 q_3^2}{12 \left(q_1^2+q_2^2+q_3^2\right)^{5/2}}.
\end{eqnarray}
In this setting $\mathcal H$ leads to a perturbed Keplerian system. Next, we express the Hamiltonian in Delaunay variables \cite{Deprit1981}. For a compact notation, we consider the eccentricity $e=\sqrt{1-G^2/L^2},$ $\eta=\sqrt{1-e^2},$ $\cos x=c_x$ and $\sin x=s_x.$
$$\bar{\mathcal H}(\ell,g,h,L,G,H)=\mathcal H_0^0+\varepsilon{\mathcal H}_1^0+\frac{\varepsilon^2}{2}{\mathcal H}_2^0+O(\varepsilon^3),$$
where,
\begin{eqnarray}
\nonumber {\mathcal H}_0^0&=&-\dfrac{1}{2L^2},\\
\nonumber {\mathcal H}_1^0&=& \frac{1}{6 G L^4 \left(e c_E -1\right)^3}\left[c_g \left(G \left(c_E-e\right) \left(3 c_h+\sqrt{3} s_h\right)+\eta  H s_E \left(\sqrt{3} c_h-3 s_h\right)\right)\right.\\
\nonumber &&\qquad\left.-s_g \left(H \left(e-c_E\right) \left(\sqrt{3} c_h-3 s_h\right)+\eta  G s_E \left(3 c_h+\sqrt{3} s_h\right)\right)\right],
\\
\nonumber{\mathcal H}_2^0&=&\frac{\left(c_E-e\right){}^2 c_g^2 H^2}{2 G^2 L^6 \left(1-e c_E\right){}^5}+\frac{c_h^2 H^2}{8 G^2 L^6 \left(1-e c_E\right){}^3}+\frac{\left(1-e^2\right) c_g^2 c_h^2 s_E^2 H^2}{8 G^2 L^6 \left(1-e c_E\right){}^5}+\dots
\end{eqnarray}
and the Delaunay variables are restricted to the domain $M=\mathbb T^3\times \Omega$ with $\Omega=\{(L,G,H):\vert H\vert <G<L,G\neq 0\}.$ The complete expression for ${\mathcal H}_2^0$ can be found in Appendix C. We are now ready to normalize using the Lie-Deprit canonical transformation
\cite{Deprit82,de},
$\mathcal L:(\ell,g,h,L,G,H)\to(\ell',g',h',L',G',H').$
This process relies on the following recurrence
\begin{equation}\label{eq:algoritmoNorm}
\mathcal H_j^i=\mathcal H_{j+1}^{i-1}+\sum_{k=0}^j\left(
\begin{array}{c}
j\\
k\\
\end{array}
\right)\{\mathcal H_{j-k}^{i-1},\mathcal W_{k+1}\},
\end{equation}
where $\bar{\mathcal{H}}=\sum_{i=0}^\infty\left(\varepsilon^i/i!\right)\mathcal H_i^0$ and $\tilde{\mathcal H}=\sum_{i=0}^\infty\left(\varepsilon^i/i!\right)\mathcal H_0^i$ are the Hamiltonian functions in the old and new variables respectively and $\mathcal W=\sum_{i=0}^\infty\left(\varepsilon^i/i!\right)\mathcal W_{i+1}$ is the generating function. 

Since we have $\mathcal{H}_0^1 = 0,$ we proceed to compute $\mathcal{H}_0^2$. The normalized Hamiltonian, up to the second order, then takes the final form:
\begin{equation}\label{HamPromediado}
\tilde{\mathcal H}=-\frac{1}{2 L'^2}-\frac{\varepsilon^2}{2}\frac{8 H'^2+\left(G'^2-H'^2\right) \left(\sqrt{3} \sin 2 h'+\cos 2 h'\right)}{24 G'^5 L'^3}+O(\varepsilon^3).
\end{equation}
In what follows, we omit the primes, as the variables will be understood from the context. We also include $\mathcal{W}_1$ for the benefit of readers who wish to follow the details for computing $\mathcal H_0^2$. 
\begin{eqnarray}
\nonumber \mathcal W_1&=&\frac{1}{6 e G L(e c_E-1)}
\left[G \left(\sqrt{3} s_h+3 c_h\right) (e s_E c_g+\eta  s_g)-H \left(\sqrt{3} c_h-3s_h\right) (\eta  c_g-e s_E s_g)\right].
\end{eqnarray}

\section{First reduced space}\label{sec:Reducido1}
The open subset of the Keplerian phase space determined by $\mathcal H_0<0,$ $G\neq0$ is a $S^1$-fiber bundle. Moreover, a coordinate system for the base space called as the reduced space is 
\begin{equation}
x=\mathbf G+L\mathbf A,\qquad y=\mathbf G-L\mathbf A,
\end{equation}
where $\mathbf G$ is the angular momentum and $\mathbf A$ is the Laplace-Runge-Lenz vector. The explicit expressions in Delaunay variables of $ x$ and $ y$ are found in \cite{Coffey1986}. Moreover, we have the following Poisson bracket relations
$$\{x_i,x_j\}=2\varepsilon_{ijk}\, x_k,\:\:\{y_i,y_j\}=2\varepsilon_{ijk}\, y_k, \:\:\{x_i,y_j\}=0,\:\:i,j,k\in\{1,2,3\},$$
{\color{black}being  $\varepsilon_{ijk}$ is the Levi-Civita symbol. With the Casimirs
$$C_1:=x_1^2+x_2^2+x_3^2=L^2, \quad  C_2:=y_1^2+y_2^2+y_3^3=L^2.$$
Therefore, considering the aforementioned constraints, the reduced space is given by $S^2 \times S^2.$ 
After some algebraic manipulations, the reduced Hamiltonian is given by
\begin{equation}\label{eq:RedHam}
\mathcal H_\ell=-\frac{1}{2 L^2} + \frac{\varepsilon^2}{2} \frac{\alpha _{x_1}^2+2 \sqrt{3} \alpha _{x_2} \alpha _{x_1}-\alpha _{x_2}^2-8 \alpha _{x_3}^2}{3 (2\alpha) ^{5/2} L^3}+O(\varepsilon^3),
\end{equation}
where $\alpha=\frac{1}{2}\Vert x+ y\Vert ^2$ and  $\partial \alpha /\partial x_i=\alpha_{x_i} $.  The equations of motion are given by
\begin{equation}\label{eq:EcuMovReduce}
\begin{split}
\dot x_1=&\frac{\varepsilon^2}{24 \sqrt{2} \alpha ^{7/2} L^3}\left[x_2 \alpha _{x_3} \left(32 \alpha +5 \alpha _{x_1}^2+10 \sqrt{3} \alpha _{x_2} \alpha _{x_1}-5 \alpha _{x_2}^2-40 \alpha _{x_3}^2\right)\right.\\
&\quad\left.+x_3 \left(-5 \alpha _{x_2} \alpha _{x_1}^2+2 \sqrt{3} \left(2 \alpha -5 \alpha _{x_2}^2\right) \alpha _{x_1}+\alpha _{x_2} \left(-4 \alpha +5 \alpha _{x_2}^2+40 \alpha _{x_3}^2\right)\right)  \right]+O(\varepsilon^3),\\
\dot x_2=&\frac{\varepsilon^2}{24 \sqrt{2} \alpha ^{7/2} L^3}\left[ 5 x_3 \alpha _{x_1}^3+5 \left(2 \sqrt{3} x_3 \alpha _{x_2}-x_1 \alpha _{x_3}\right) \alpha _{x_1}^2-\left(5 x_3 \alpha _{x_2}^2+10 \sqrt{3} x_1 \alpha _{x_3} \alpha _{x_2}\right.\right.\\
&\quad\left.\left.+4 x_3 \left(\alpha +10 \alpha _{x_3}^2\right)\right) \alpha _{x_1}-4 \sqrt{3} \alpha  x_3 \alpha _{x_2}+5 x_1 \alpha _{x_2}^2 \alpha _{x_3}+8 x_1 \alpha _{x_3} \left(5 \alpha _{x_3}^2-4 \alpha \right) \right]+O(\varepsilon^3),\\
\dot x_3=&\frac{\varepsilon^2}{24 \sqrt{2} \alpha ^{7/2} L^3}\left[  -5 x_2 \alpha _{x_1}^3+5 \left(x_1-2 \sqrt{3} x_2\right) \alpha _{x_2} \alpha _{x_1}^2+\left(5 \left(2 \sqrt{3} x_1+x_2\right) \alpha _{x_2}^2+40 x_2 \alpha _{x_3}^2\right.\right.\\
&\quad\left.\left.-4 \sqrt{3} \alpha  x_1+4 \alpha  x_2\right) \alpha _{x_1}+\alpha _{x_2} \left(4 \alpha  \left(x_1+\sqrt{3} x_2\right)-5 x_1 \left(\alpha _{x_2}^2+8 \alpha _{x_3}^2\right)\right) \right]+O(\varepsilon^3),\\
\dot y_1=&\frac{\varepsilon^2}{24 \sqrt{2} \alpha ^{7/2} L^3}\left[-5 \left(\alpha _{x_1}^2+2 \sqrt{3} \alpha _{x_2} \alpha _{x_1}-\alpha _{x_2}^2-8 \alpha _{x_3}^2\right)\left(y_3 \alpha _{y_2}-y_2 \alpha _{y_3}\right) \right.\\
&\quad\left.+4 \alpha  \left(\sqrt{3} y_3 \alpha _{x_1}-y_3 \alpha _{x_2}+8 y_2 \alpha _{x_3}\right) \right]+O(\varepsilon^3),\\
\dot y_2=&\frac{\varepsilon^2}{24 \sqrt{2} \alpha ^{7/2} L^3}\left[  5 \left(\alpha _{x_1}^2+2 \sqrt{3} \alpha _{x_2} \alpha _{x_1}-\alpha _{x_2}^2-8 \alpha _{x_3}^2\right) \left(y_3 \alpha _{y_1}-y_1 \alpha _{y_3}\right)\right.\\
&\quad\left.-4 \alpha  \left(y_3 \alpha _{x_1}+\sqrt{3} y_3 \alpha _{x_2}+8 y_1 \alpha _{x_3}\right) \right]+O(\varepsilon^3),\\
\dot y_3=&\frac{\varepsilon^2}{24 \sqrt{2} \alpha ^{7/2} L^3}\left[ -5 \left(\alpha _{x_1}^2+2 \sqrt{3} \alpha _{x_2} \alpha _{x_1}-\alpha _{x_2}^2-8 \alpha _{x_3}^2\right) \left(y_2 \alpha _{y_1}-y_1 \alpha _{y_2}\right)\right.\\
&\quad\left.+4 \alpha  \left(\left(y_2-\sqrt{3} y_1\right) \alpha _{x_1}+\left(y_1+\sqrt{3} y_2\right) \alpha _{x_2}\right)  \right]+O(\varepsilon^3).
\end{split}
\end{equation}

\subsection{Relative Equilibria and Stability}
The system \eqref{eq:EcuMovReduce} has several families of equilibria; however, we focus on the isolated equilibria $\mathbf E_i=(x^i,y^i)$, given by
\begin{equation}\label{eq:rev1eq}
\mathbf E_{i} =\delta_i(0,0,1,0,0,1),\quad \mathbf E_{i+2} =\delta_i \left(\dfrac{\sqrt{3}}{2},\dfrac{1}{2},0,\dfrac{\sqrt{3}}{2},\dfrac{1}{2},0\right),\quad \mathbf E_{i+4} = \delta_i\left(\dfrac{1}{2},-\dfrac{\sqrt{3}}{2},0,\dfrac{1}{2},-\dfrac{\sqrt{3}}{2},0\right),
\end{equation}
where $i=1,2$ and $\delta_i=(-L)^{i+1}$. All the equilibria $\mathbf E_i$ correspond to circular orbits. 
Furthermore, the equilibria $\mathbf E_{1,2}$ correspond to direct and retrograde circular equatorial orbits of the unperturbed system, while the remaining equilibria $\mathbf E_{3,4}$ and $\mathbf E_{5,6}$ correspond to polar orbits. The reconstruction of these equilibria in the complete system is depicted in Figure~\ref{fig:2rec}.

\begin{Theorem}\label{teo1}
The relative equilibira $\mathbf E_i$ reconstruct into circular periodic orbits in the full system, with period $T=2\pi L^3+O(\varepsilon^3).$
\end{Theorem}
\begin{proof} Consider the relative equilibria $\mathbf E_1$ and $\mathbf E_2,$ and introduce the following local transformation
\begin{equation}\label{eq:XtoQ}
\begin{split}
Q_1=\dfrac{x_2}{\sqrt{L\pm  x_3}},\quad Q_2=\dfrac{ y_2}{\sqrt{L\pm y_3}},\quad
P_1=\mp\dfrac{ x_1}{\sqrt{L\pm  x_3}},\quad P_2=\mp\dfrac{ y_1}{\sqrt{L\pm y_3}}.
\end{split}
\end{equation}
The inverse transformation is given by
\begin{equation}\nonumber
\begin{split}
 x_1=&\mp P_1\sqrt{2L-P_1^2-Q_1^2},\qquad  x_2= Q_1\sqrt{2L-P_1^2-Q_1^2},\qquad  x_3=\pm L\mp (P_1^2+Q_1^2),\\
y_1=&\mp P_2\sqrt{2L-P_2^2-Q_2^2},\qquad  y_2= Q_2\sqrt{2L-P_2^2-Q_2^2},\qquad  y_3=\pm L\mp (P_2^2+Q_2^2).
\end{split}
\end{equation}
Note that the upper sign corresponds to the first equilibrium, and the lower sign to the second. The variables $(Q_1, Q_2, P_1, P_2)$ are canonical, where $Q_1$ and $Q_2$ represent the coordinates, and $P_1$ and $P_2$ are their corresponding momenta. The equilibria in the new variables become $(0,0,0,0).$
The resulting Hamiltonian is obtained by expressing $\mathcal H_\ell$, given in \eqref{eq:RedHam} in terms of $(Q_1,Q_2,P_1,P_2),$ dropping constant terms and, rescaling the time by $d\tau=\varepsilon^2/2 dt.$ 
\begin{equation}\label{Hbarcir12}
\mathcal H_\ell^{1,2}(Q_1,Q_2,P_1,P_2)=\frac{-\mathfrak a^2\mp2 \sqrt{3} \mathfrak a \mathfrak b+\mathfrak b^2-8 \left(-2 L+v^2+w^2\right)^2}{3 L^3 \left(\mathfrak a^2+\mathfrak b^2+\left(-2 L+v^2+w^2\right)^2\right)^{5/2}}+O(\varepsilon),
\end{equation}
where, $ v=(Q_1,P_1),$ $ w=(Q_2,P_2),$ $R_{ v}=\sqrt{2L- v^2},$ $R_{ w}=\sqrt{2L- w^2},$ $ \mathfrak a=Q_1R_ v+Q_2R_ w$ and $\mathfrak b=P_1R_ v+P_2R_ w.$ The Hessian matrix of the vector field, evaluated at each equilibrium point, is: $$H_{1,2}=\left[
\begin{array}{cccc}
 -\frac{5}{24 L^7} & \frac{19}{24 L^7} & \mp\frac{1}{8 \sqrt{3} L^7} & \mp\frac{1}{8 \sqrt{3} L^7} \\
 \frac{19}{24 L^7} & -\frac{5}{24 L^7} & \mp\frac{1}{8 \sqrt{3} L^7} & \mp\frac{1}{8 \sqrt{3} L^7} \\
 \mp\frac{1}{8 \sqrt{3} L^7} & \mp\frac{1}{8 \sqrt{3} L^7} & -\frac{1}{8 L^7} & \frac{7}{8 L^7} \\
 \mp\frac{1}{8 \sqrt{3} L^7} & \mp\frac{1}{8 \sqrt{3} L^7} & \frac{7}{8 L^7} & -\frac{1}{8 L^7} \\
\end{array}
\right],$$
and $det(H_{1,2})=5/(12 L^{28}).$ 
Thus, the nondegeneracy condition of Reeb’s theorem \cite{Palacian2008} is satisfied for each critical point, leading to periodic solutions in the original system. Moreover, the characteristic multiplies of these orbits are
$$1,1,1+\varepsilon^2 \dfrac{i}{ L^7}T+O(\varepsilon^3),1+\varepsilon^2 \dfrac{-i}{ L^7}T+O(\varepsilon^3),1+\varepsilon^2\sqrt{\frac{5}{3}} \dfrac{i}{2 L^7}T+O(\varepsilon^3),1+\varepsilon^2 \sqrt{\frac{5}{3}}\dfrac{-i}{2 L^7}T+O(\varepsilon^3).$$

The analyses for the equilibria $\mathbf E_i$ for $i=3,4,5,6$ are analogous. Precisely, the Hessian matrix is the same for each of the following pairs of equilibria $(\mathbf E_{3}, \mathbf E_4)$ and $(\mathbf E_{5}, \mathbf E_6)$, we denote them by $H_{3,4}$ and $H_{56}$, which are respectively given by
$$\left[
\begin{array}{cccc}
 -\frac{31}{192 L^7} & -\frac{73}{192 L^7} &\pm \frac{9 \sqrt{3}}{64 L^7} & \pm\frac{15 \sqrt{3}}{64 L^7} \\[1ex]
 -\frac{73}{192 L^7} & -\frac{31}{192 L^7} &\pm \frac{15 \sqrt{3}}{64 L^7} &\pm \frac{9 \sqrt{3}}{64 L^7} \\[1ex]
\pm \frac{9 \sqrt{3}}{64 L^7} &\pm \frac{15 \sqrt{3}}{64 L^7} & -\frac{85}{192 L^7} & -\frac{163}{192 L^7} \\[1ex]
 \pm\frac{15 \sqrt{3}}{64 L^7} &\pm \frac{9 \sqrt{3}}{64 L^7} & -\frac{163}{192 L^7} & -\frac{85}{192 L^7} \\[1ex]
\end{array}
\right],
\quad \left[
\begin{array}{cccc}
 -\frac{107}{192 L^7} & -\frac{29}{192 L^7} & \pm\frac{37}{64 \sqrt{3} L^7} &\pm \frac{19}{64 \sqrt{3} L^7} \\[1ex]
 -\frac{29}{192 L^7} & -\frac{107}{192 L^7} & \pm\frac{19}{64 \sqrt{3} L^7} &\pm \frac{37}{64 \sqrt{3} L^7} \\[1ex]
 \pm\frac{37}{64 \sqrt{3} L^7} & \pm\frac{19}{64 \sqrt{3} L^7} & -\frac{11}{64 L^7} & \frac{3}{64 L^7} \\[1ex]
 \pm\frac{19}{64 \sqrt{3} L^7} & \pm\frac{37}{64 \sqrt{3} L^7} & \frac{3}{64 L^7} & -\frac{11}{64 L^7} \\[1ex]
\end{array}
\right].$$
Hence, we have $ det(H_{3,4})=-5/3\,det(H_{5,6})=\frac{5}{288 L^{28}}\neq0.$
\end{proof}

\begin{figure}[h]
    \centering
    \subfigure[]{ \includegraphics[width=147pt]{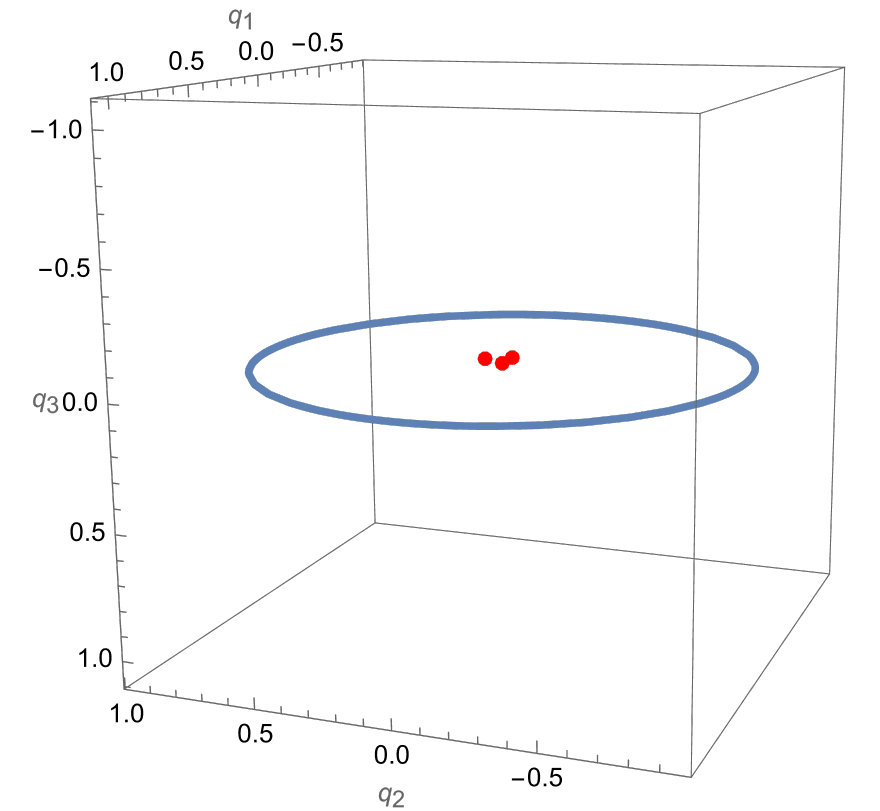}} 
    \subfigure[]{ \includegraphics[width=147pt]{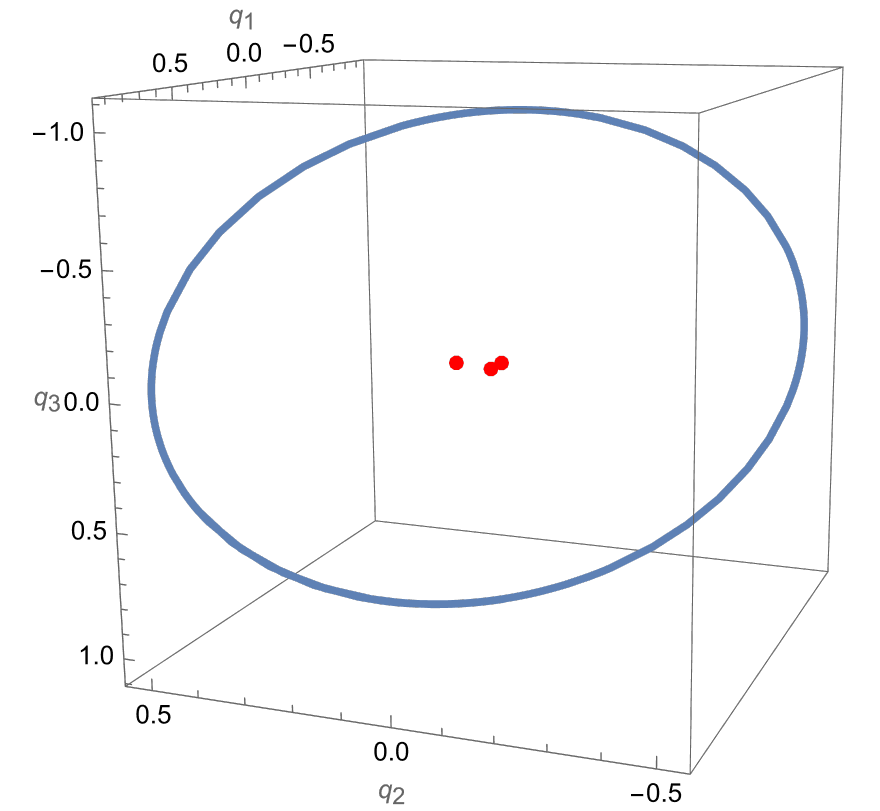}} 
    \subfigure[]{ \includegraphics[width=147pt]{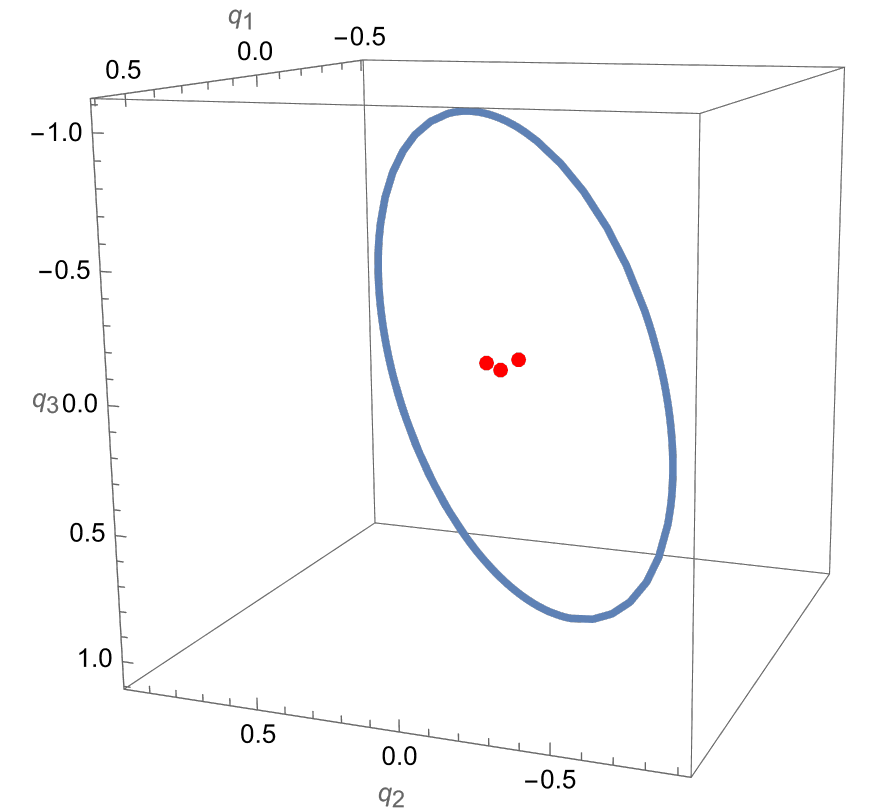}}
    \caption{ \small Simulation in the configuration space of the full system \eqref{eq:Ham3centros} of the periodic orbits reconstructed from \eqref{eq:rev1eq}. We consider $L=1$ and $\varepsilon=0.1$. $\mathbf{E}_{1,2}$ correspond with equatorial orbits illustrated in (a), while $\mathbf{E}_{3,4}$ and $\mathbf{E}_{5,6}$ reconstruct to the polar orbits showed in (b) and (c) respectively.}
    \label{fig:2rec}
\end{figure}

\begin{Theorem}\label{Teo_estabilidad_1}
The periodic solutions reconstructed from the relative equilibria $\mathbf{E}_1,$ $\mathbf E_2,$ $\mathbf E_3$ and $\mathbf E_4$ are parametrically stable. In contrast, the periodic solutions generated by $\mathbf E_5$ and $\mathbf E_6$ are unstable.\end{Theorem}
\begin{proof}

Consider the relative equilibria $\mathbf E_{1,2}$. The linearization matrix and its characteristic polynomial are given by
$$A_{1,2}=\left[
\begin{array}{cccc}
 \mp\frac{1}{8 \sqrt{3} L^7} & \mp\frac{1}{8 \sqrt{3} L^7} & -\frac{1}{8 L^7} & \frac{7}{8 L^7} \\[1ex]
 \mp\frac{1}{8 \sqrt{3} L^7} & \mp\frac{1}{8 \sqrt{3} L^7} & \frac{7}{8 L^7} & -\frac{1}{8 L^7} \\[1ex]
 \frac{5}{24 L^7} & -\frac{19}{24 L^7} &\pm \frac{1}{8 \sqrt{3} L^7} &\pm \frac{1}{8 \sqrt{3} L^7} \\[1ex]
 -\frac{19}{24 L^7} & \frac{5}{24 L^7} &\pm \frac{1}{8 \sqrt{3} L^7} &\pm \frac{1}{8 \sqrt{3} L^7} \\[1ex]
\end{array}
\right],\quad C(\lambda)=\lambda ^4+\frac{5}{12 L^{28}}+\frac{17 \lambda ^2}{12 L^{14}}.$$
The matrix $A_{1,2}$ is diagonalizable, with eigenvalues $\pm \alpha_1i$ and $\pm \alpha_2i,$ where
$\alpha_1=1/L^7,$ $\alpha_2=\sqrt{5}/(2\sqrt{3}L^7).$
The quadratic Hamiltonian associated with $A_{1,2}$ in some symplectic coordinates $(q_1,q_2,p_1,p_2)$ is given by
\begin{equation}\nonumber
\mathcal K_2^{1,2}=\frac{-6(q_1^2+p_1^2)+\sqrt{15}(q_2^2+ p_2^2)}{6 L^7}.
\end{equation}
The maximal real linear subspace associated to $\mathcal \alpha_1$ and $\mathcal \alpha_2$ is $V_1=Ker(A_{1,2}-\alpha_1iId)$ and $V_2=Ker(A_{1,2}-\alpha_2iId)$, where
$$V_1=\left<\{0, 0, -1, 1\},\{-1, 1, 0, 0\}\right>,\ V_2=\left<\left\{\pm\sqrt{3}/7,\pm\sqrt{3}/7,1,1\right\},\left\{-2\sqrt{15}/7,-2\sqrt{15}/7,0,0\right\}\right>.$$
Therefore, $\mathcal K^{1,2}_2|_{V_1}$ and $\mathcal K_2^{1,2}|_{V_2}$ have a definite sign. As a consequence, $A_{1,2}$ is parametrically stable. By the Theorem~\ref{Teo:inestable}, the direct and retrograde equatorial periodic solution are stable.

Now, consider the relative equilibria $\mathbf E_{3,4}.$ The matrix of linearization around them and its corresponding characteristic polynomial take the form
$$A_{3,4}=\left[
\begin{array}{cccc}
 \frac{\pm9 \sqrt{3}}{64 L^7} &\pm \frac{15 \sqrt{3}}{64 L^7} & -\frac{85}{192 L^7} & -\frac{163}{192 L^7} \\[1ex]
 \frac{\pm15 \sqrt{3}}{64 L^7} &\pm \frac{9 \sqrt{3}}{64 L^7} & -\frac{163}{192 L^7} & -\frac{85}{192 L^7} \\[1ex]
 \frac{31}{192 L^7} & \frac{73}{192 L^7} & \mp\frac{9 \sqrt{3}}{64 L^7} & \mp\frac{15 \sqrt{3}}{64 L^7} \\[1ex]
 \frac{73}{192 L^7} & \frac{31}{192 L^7} & \mp\frac{15 \sqrt{3}}{64 L^7} & \mp\frac{9 \sqrt{3}}{64 L^7} \\[1ex]
\end{array}
\right],\quad C(\lambda)=\lambda ^4+\frac{5}{288 L^{28}}+\frac{49 \lambda ^2}{144 L^{14}}.$$
Therefore, $A_{3,4}$ is diagonalizable, with eigenvalues $\pm\alpha_1 i$ and $\pm\alpha_2i,$ where
$\alpha_1=\sqrt{5/2}/(3 L^7)$ and $ \alpha_2=1/(4 L^7).$
The quadratic Hamiltonian associated with $A_{3,4}$ in some symplectic coordinates $(q,p)$ is given by
\begin{equation}\nonumber
\mathcal K_2^{3,4}=\frac{-2 \sqrt{10} (q_1^2+p_1^2)+3 (q_1^2+p_2^2)}{24 L^7}.
\end{equation}
The maximal real linear subspace associated to $\mathcal \alpha_1$ and $\mathcal \alpha_2$ is $V_3=Ker(A_{3,4}-\alpha_1iId),$ $V_4=Ker(A_{3,4}-\alpha_2iId)$ respectively, where
\begin{equation}\nonumber
\begin{split}
V_3&=\left<\left\{\pm9 \sqrt{3}/13,\pm9 \sqrt{3}/13,1,1\right\},\left\{4 \sqrt{10}/13,4 \sqrt{10}/13,0,0\right\}\right>,\\
V_4&=\left<\left\{\mp3 \sqrt{3}/7,\pm3 \sqrt{3}/7 ,-1,1\right\},\left\{8/7,-8/7,0,0\right\}\right>,
\end{split}
\end{equation}
and $\mathcal K_2^{3,4}|_{V_3},$ $\mathcal K_2^{3,4}|_{V_4},$ have a definite sign. $A_{3,4}$ is parametrically stable and by the Theorem~\ref{Teo:inestable} the direct and retrograde polar periodic solution are stable.

Finally, if we consider the equilibria $\mathbf E_{5,6}.$ The linearization matrix and its characteristic polynomial are given by
$$A_{5,6}=\left[
\begin{array}{cccc}
 \pm\frac{37}{64 \sqrt{3} L^7} & \pm\frac{19}{64 \sqrt{3} L^7} & -\frac{11}{64 L^7} & \frac{3}{64 L^7} \\[1ex]
 \pm\frac{19}{64 \sqrt{3} L^7} &\pm \frac{37}{64 \sqrt{3} L^7} & \frac{3}{64 L^7} & -\frac{11}{64 L^7} \\[1ex]
 \frac{107}{192 L^7} & \frac{29}{192 L^7} & \mp\frac{37}{64 \sqrt{3} L^7} & \mp\frac{19}{64 \sqrt{3} L^7} \\[1ex]
 \frac{29}{192 L^7} & \frac{107}{192 L^7} & \mp\frac{19}{64 \sqrt{3} L^7} & \mp\frac{37}{64 \sqrt{3} L^7} \\[1ex]
\end{array}
\right],\quad C(\lambda)=\lambda ^4-\frac{1}{96 L^{28}}-\frac{5 \lambda ^2}{48 L^{14}}$$
and has one real eigenvalue. Consequently, by the Theorem~\ref{Teo:inestable}, the polar periodic solution generated by $\mathbf E_{5,6}$ are unstable in the Lyapunov sense.
\end{proof}

\begin{figure}[h]
    \centering
    \subfigure[]{ \includegraphics[width=147pt]{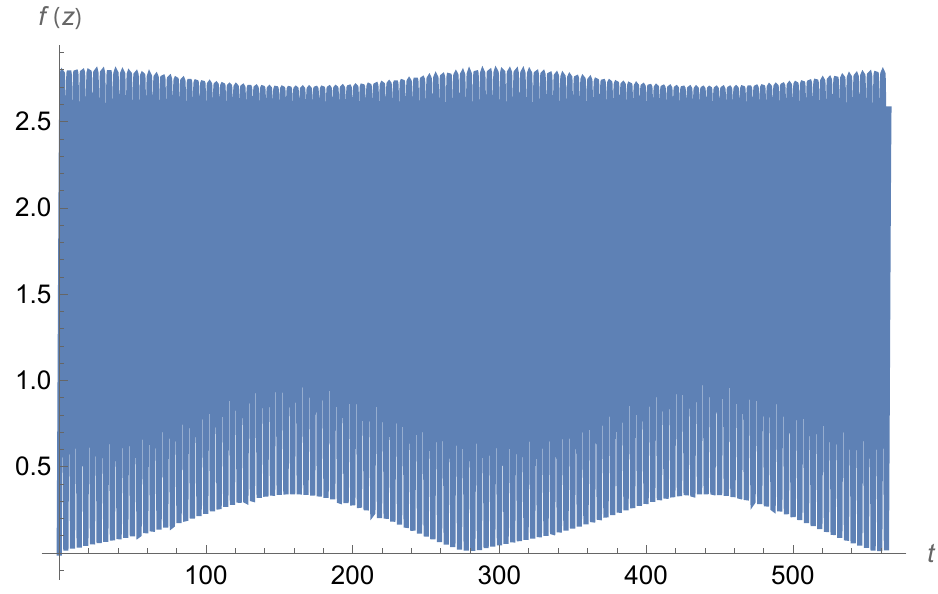}} 
    \subfigure[]{ \includegraphics[width=147pt]{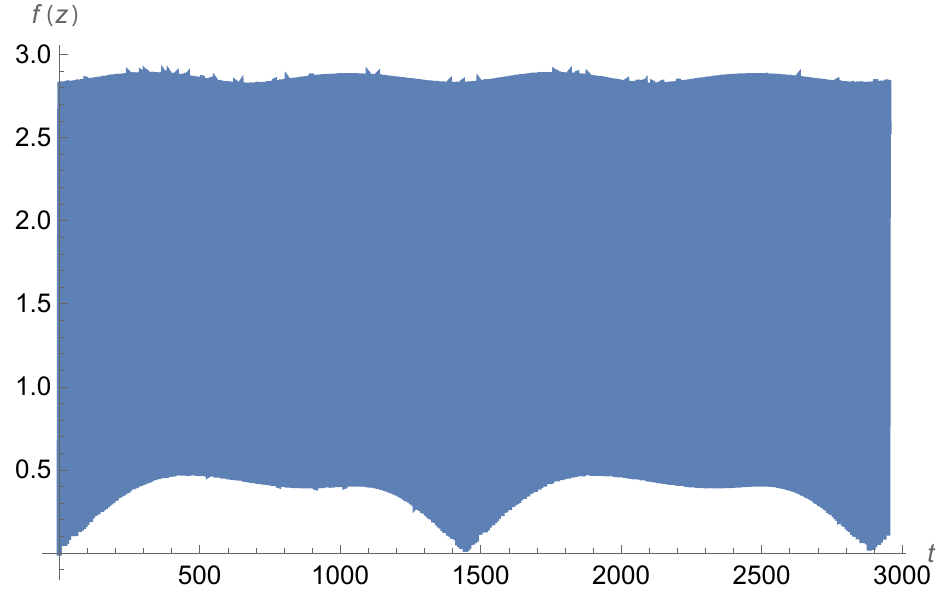}} 
    \subfigure[]{ \includegraphics[width=147pt]{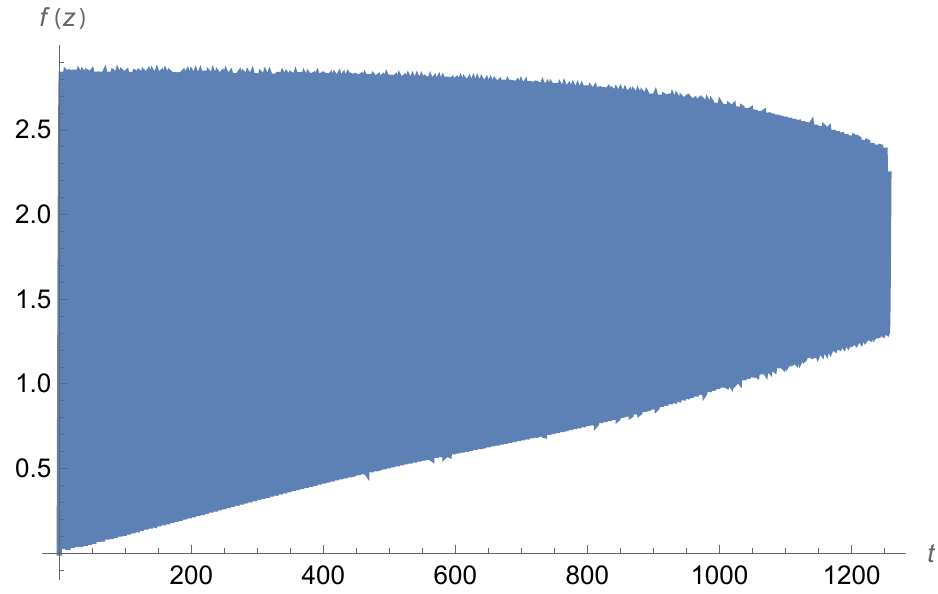}}
    \caption{\small Plot of the function $f(z) = \Vert z(t_0) - z(t) \Vert$, with $z=(q,p)$. }
    \label{fig:3rec}
\end{figure}

The Figure~\ref{fig:3rec} provides insight into how closely the solution approximates a periodic solution within the full system \eqref{eq:Ham3centros}. Specifically, the function plotted shows the distance of the solution \( (q(t),p(t)) \), from the initial condition \( (q(0),p(0)) \) at each point in time. Additionally, it reveals the stability of equilibria \( \mathbf{E}_{1,2} \) and \( \mathbf{E}_{3,4} \), as demonstrated in Theorem~\ref{Teo_estabilidad_1}, along with the instability exhibited by the equilibrium points \( \mathbf{E}_{5,6} \). This visual representation allows us to observe these stability characteristics directly, reinforcing the theoretical results.

\subsection{KAM 3-Tori}
In this section, we demonstrate the existence of invariant three-dimensional tori surrounding the near-circular periodic solutions $\mathbf{E}_{1,2}$ and $\mathbf{E}_{3,4}$. A key element is the use of a KAM theorem by Han, Li, and Yi \cite{hanliyi}, which addresses Hamiltonians with high-order proper degeneracy, where perturbations occur at different orders.

Under the assumptions of the Theorem~\ref{Teo_estabilidad_1} we obtain the following result.
\begin{Theorem}
There are families of invariant KAM 3-tori of the full system around the near circular periodic solutions $\mathbf E_{1,2}$ and $\mathbf E_{3,4}$ given in Theorem~\ref{teo1}. These invariant tori form a majority in the sense that the measure of the complement of their union is of the order $O(\varepsilon^{\delta/4}),$ with $1<\delta<\frac{1}{5}.$
\end{Theorem}
\begin{proof}
Consider the Hamiltonian $\mathcal H_\ell$ given in \eqref{eq:RedHam}, and apply the symplectic local transformation \eqref{eq:XtoQ}, arriving at $\mathcal H_\ell^{1,2}$ as given in \eqref{Hbarcir12}. The equilibria $\mathbf E_{1,2}$ in the new variables become $(0,0,0,0).$ Then, applying the following $\varepsilon^{-1/4}$-symplectic transformation on $\mathcal H_\ell^{1,2},$ we obtain
\begin{equation}
\begin{split}
&q_1=\frac{\varepsilon^{1/8}}{70}\left[35 \sqrt{2} P_1\pm\sqrt[4]{15} \sqrt{7} \left(\sqrt{5} Q_2-10 P_2\right)\right],\quad p_1=\varepsilon^{1/8}\left(\frac{\sqrt{7} Q_2}{2 \sqrt[4]{15}}-\frac{Q_1}{\sqrt{2}}\right),\\
&q_2=\frac{\varepsilon^{1/8}}{70}\left[-35 \sqrt{2} P_1\pm\sqrt[4]{15} \sqrt{7} \left(\sqrt{5} Q_2-10 P_2\right)\right],\quad p_2=\varepsilon^{1/8}\left(\frac{\sqrt{7} Q_2}{2 \sqrt[4]{15}}+\frac{Q_1}{\sqrt{2}}\right),
\end{split}
\end{equation}
we obtain $
\mathcal K^{1,2}(q_1,q_2,p_1,p_2)=\varepsilon^{-1/4}\mathcal K_1^{1,2}+\mathcal K_2^{1,2}+\varepsilon^{1/4}\mathcal K_4^{1,2}+O(\varepsilon ^{1/2}),
$ where,
\begin{eqnarray}
\nonumber \mathcal K_1^{1,2}&=&-\frac{1}{3 L^6},\\
\nonumber \mathcal K_2^{1,2}&=&-\frac{6(q_1^2-p_1^2)+\sqrt{15}(q_2^2+ p_2^2)}{12 L^7},\\
\nonumber \mathcal K_4^{1,2}&=& -\frac{1}{1680 L^8}\left[\mp p_1^2 \left(95 \sqrt{3} p_2 q_2\pm185 \sqrt{15} p_2^2\mp1680 q_1^2\pm361 \sqrt{15} q_2^2\right)\right.\\
\nonumber &&\quad+2 p_1 q_1 \left(685 p_2 q_2\pm15 \sqrt{5} p_2^2\mp57 \sqrt{5} q_2^2\right)+10 \left(p_2^2+q_2^2\right) \left(\mp3 \sqrt{5} p_2 q_2+15 p_2^2+13 q_2^2\right)\\
\nonumber &&\quad\left.+\sqrt{3} q_1^2 \left(95 p_2 q_2\mp375 \sqrt{5} p_2^2\mp199 \sqrt{5} q_2^2\right)+840 p_1^4+840 q_1^4\right].
\end{eqnarray}
Next, we normalize the quartic terms $\mathcal K_4^{1,2}$ using a complex symplectic transformation. Applying the Lie-Deprit method and returning to Cartesian coordinates, we obtain $\bar{\mathcal K}_1^{1,2}= \mathcal K_1^{1,2},$ $\bar{\mathcal K}_2^{1,2}={\mathcal K}_2^{1,2}$ and
\begin{eqnarray}
\nonumber \bar{\mathcal K}_4^{1,2}&=&-\frac{\left(p_1^2+q_1^2\right){}^2}{2 L^8}+\frac{\sqrt{\frac{5}{3}} \left(p_2^2+q_2^2\right) \left(p_1^2+q_1^2\right)}{2 L^8}-\frac{\left(p_2^2+q_2^2\right){}^2}{12 L^8}.\end{eqnarray}
As the next step, we introduce action-angle variables $(I_1,I_2,\phi_1,\phi_2),$ which are defined by the transformation
\begin{equation}\label{accionangulo}
\begin{split}
&q_1 = \sqrt{2 I_1} \cos\phi_1, \quad q_2 = \sqrt{2 I_2} \cos\phi_2,\quad
p_1 = \sqrt{2 I_1} \sin\phi_1, \quad p_2 = \sqrt{2 I_2} \sin\phi_2,
\end{split}
\end{equation}
incorporate the terms associated to the action $L$ dropped in the previous steps and undo the time scalings, arriving at 
\begin{equation}\nonumber
\tilde{\mathcal K}^{1,2}=-\frac{1}{2 L^2}-\frac{\eta^8}{3 L^6}-\frac{\eta^9(6 I_1-\sqrt{15}I_2)}{6 L^7}-\frac{\eta^{10}(6 I_1^2-2 \sqrt{15} I_1 I_2+I_2^2)}{3 L^8}+O(\eta^{11}),
\end{equation}
where $\eta=\varepsilon^{1/4}.$ Finally, we are in the position to apply Han, Li and Yi's Theorem~\ref{Hanliyi} (see Appendix~\ref{apB}) taking $a=3,$ $n=3,$ $m_1=8,$ $m_2=9,$ $m_3=10,$ $n_0=n_1=1,$ $n_2=n_3=3,$ $I^{n_0}=I^{n_1}=L,$ $I^{n_2}=I^{n_3}=(L,I_1,I_2),$ $\bar I^{n_0}=\bar I^{n_1}=L,$ $\bar I^{n_2}=\bar I^{n_3}=(I_1,I_2).$ In this case, the vector of frequencies has a dimension 6 and is given by
\begin{equation}
\begin{split}
\Omega(I)&=
\left(\frac{1}{L^3},\frac{2}{L^7},-\frac{1}{L^7},\frac{\sqrt{\frac{5}{3}}}{2 L^7},-\frac{2 \left(6 I_1-\sqrt{15} I_2\right)}{3 L^8},\frac{2 \left(\sqrt{15} I_1-I_2\right)}{3 L^8}\right).
\end{split}
\end{equation}
The corresponding matrix with dimension $6\times4$ whose columns are $\Omega(I),$ $\partial \Omega(I)/\partial L,$ $\partial\Omega(I)/\partial I_1$ and $\partial\Omega(I)/\partial I_2$ reads
\begin{equation}
M_{\Omega}(L,I_1,I_2)=\left[
\begin{array}{cccc}
 \frac{1}{L^3} & -\frac{3}{L^4} & 0 & 0 \\[1ex]
 \frac{2}{L^7} & -\frac{14}{L^8} & 0 & 0 \\[1ex]
 -\frac{1}{L^7} & \frac{7}{L^8} & 0 & 0 \\[1ex]
 \frac{\sqrt{5/3}}{2 L^7} & -\frac{7 \sqrt{5/3}}{2 L^8} & 0 & 0 \\[1ex]
 -\frac{2 \left(6 I_1-\sqrt{15} I_2\right)}{3 L^8} & \frac{16 \left(6 I_1-\sqrt{15} I_2\right)}{3 L^9} & -\frac{4}{L^8} & \frac{2 \sqrt{5/3}}{L^8} \\[1ex]
 \frac{2 \left(\sqrt{15} I_1-I_2\right)}{3 L^8} & \frac{16 \left(I_2-\sqrt{15} I_1\right)}{3 L^9} & \frac{2 \sqrt{5/3}}{L^8} & -\frac{2}{3 L^8} \\[1ex]
\end{array}
\right]
\end{equation}
Since the rank of $M_\Omega$ is equal to four, there exist a 3-tori surrounding each periodic solution generated by the isolates equilibria $\mathbf E_{1,2}.$ 

For the case of $\mathbf E_{3,4},$ consider the Hamiltonian function $\mathcal H_\ell$ given in \eqref{eq:RedHam}, and apply the local transformation \eqref{eq:XtoQ}, arriving at $\mathcal H_\ell^{1,2}$ as given in \eqref{Hbarcir12}. The equilibria $\mathbf E_{3,4}$ in the new variables become 
$$ \mathbf E_{3,4}=\left(\pm\frac{\sqrt{L}}{2},\pm\frac{\sqrt{L}}{2},-\frac{1}{2} \sqrt{3L},-\frac{1}{2} \sqrt{3L}\right).$$
We move the system so that $\mathbf{E}_{3,4}$ is positioned at the origin. The Hamiltonian function of the translated system become
\begin{equation}\nonumber
\mathcal H_\ell^{3,4}=\frac{8 \left(-\mathfrak a^2+2 \sqrt{3} \mathfrak a \mathfrak c-32 \mathfrak b^2+\mathfrak c^2\right)}{3 L^3 \left(\mathfrak a^2+4 \mathfrak b^2+\mathfrak c^2\right)^{5/2}}+O(\varepsilon),
\end{equation}
where $ v=(Q_1,P_1)$ and $ w=(Q_2,P_2),$
\begin{eqnarray}
\nonumber \mathfrak a&=&\sqrt{L} (R_v+R_w)\pm2 (Q_1 R_v+Q_2 R_w),\quad \\
\nonumber \mathfrak b&=&\pm\sqrt{L} \left(Q_1+Q_2\mp\sqrt{3}(P_1+P_2)\right)+v^2+w^2,\\
\nonumber \mathfrak c&=&\sqrt{3L} (R_v+R_w)-2 (P_1 R_v+P_2 R_w),
\end{eqnarray}
and
\begin{equation}\nonumber
R_{ v}=\sqrt{L+\sqrt{L}(\sqrt{3}P_2\mp Q_2)- v^2},\quad R_{ w}=\sqrt{L+\sqrt{L}(\sqrt{3}P_1\mp Q_1)- w^2}.
\end{equation}
Next, applying the following $\varepsilon^{-1/4}$ symplectic transformation
\begin{equation}\nonumber
\begin{split}
&q_1=\varepsilon^{1/8}\left(\frac{2 p_2}{\sqrt{7}}\pm\frac{9 \sqrt{3/13} q_1}{2\ \sqrt[4]{2^3} \sqrt[4]{5}}-\frac{\sqrt[4]{2^3} \sqrt[4]{5} p_1}{\sqrt{13}}\mp\frac{3}{4} \sqrt{\frac{3}{7}} q_2\right)
 ,\quad p_1=\frac{\varepsilon^{1/8}}{20} \left(\sqrt[4]{2} 5^{3/4} \sqrt{13} q_1-5 \sqrt{7} q_2\right),\\
&q_2=\varepsilon^{1/8}\left(-\frac{2 p_2}{\sqrt{7}}\pm\frac{9 \sqrt{3/13} q_1}{2\ \sqrt[4]{2^3} \sqrt[4]{5}}-\frac{\sqrt[4]{2^3} \sqrt[4]{5} p_1}{\sqrt{13}}\pm\frac{3}{4} \sqrt{\frac{3}{7}} q_2\right)
 ,\quad p_2=\frac{\varepsilon^{1/8}}{20} \left(\sqrt[4]{2} 5^{3/4} \sqrt{13} q_1+5 \sqrt{7} q_2\right),\\
\end{split}
\end{equation}
on $\mathcal H_\ell^{3,4}$, we obtain
\begin{equation}
\mathcal K^{3,4}(q_1,q_2,p_1,p_2)=\varepsilon^{-1/4}\mathcal K_1^{3,4}+\mathcal K_2^{3,4}+\varepsilon^{1/8}\mathcal K_3^{3,4}+\varepsilon^{1/4}\mathcal K_4^{3,4}+O(\varepsilon ^{1/2}),
\end{equation}
where,
\begin{eqnarray}
\nonumber \mathcal K_{1}^{3,4}&=&\frac{1}{12 L^6},\\
\nonumber \mathcal K_{2}^{3,4}&=&\frac{-2 \sqrt{10} (q_1^2+p_1^2)+3 (q_1^2+p_2^2)}{24 L^7},\\
\nonumber \mathcal K_{3}^{3,4}&=&\frac{1}{2184\ 10^{3/4} \sqrt{13} L^{15/2}}\left[
\pm5740 \sqrt{5} p_1 q_1^2\pm65 \sqrt{2} p_1 \left(\mp155 \sqrt{3} p_2 q_2+76 p_2^2+274 q_2^2\right)\right.\\
\nonumber &&\quad +5810 \sqrt{6} q_1^3\left.\pm4760 \sqrt{5} p_1^3-7420 \sqrt{6} p_1^2 q_1+26 \sqrt{5} q_1 \left(\mp317 p_2 q_2+207 \sqrt{3} p_2^2-32 \sqrt{3} q_2^2\right)
\right],\\
\nonumber \mathcal K_4^{3,4}&=&\frac{1}{1135680 L^8}\left[ -5 p_1^2 \left(13 \sqrt{10} \left(\mp416 \sqrt{3} p_2 q_2+3336 p_2^2+4665 q_2^2\right)+151256 q_1^2\right)\right.\\
\nonumber &&\quad \pm63280 \sqrt{30} p_1^3 q_1+20 p_1 q_1 \left(13351 p_2 q_2\mp1404 \sqrt{3} p_2^2\pm7574 \sqrt{30} q_1^2\pm6773 \sqrt{3} q_2^2\right)\\
\nonumber &&\quad -182 \sqrt{10} q_1^2 \left(\mp78 \sqrt{3} p_2 q_2+1741 p_2^2+1323 q_2^2\right)-845 \left(\mp368 \sqrt{3} p_2^3 q_2+504 p_2^2 q_2^2\right.\\
\nonumber &&\quad \left.\left.\mp424 \sqrt{3} p_2 q_2^3+38 p_2^4+137 q_2^4\right)-78400 p_1^4-454230 q_1^4 \right].
\end{eqnarray}
Next, we normalize the system using a complex symplectic transformation, applying the Lie-Deprit method and returning to Cartesian coordinates we obtain $\bar{\mathcal K}_1^{3,4}= \mathcal K_1^{3,4},$ $\bar{\mathcal K}_2^{3,4}={\mathcal K}_2^{3,4},$ $\bar{\mathcal K}_3^{3,4}=0$ and
\begin{eqnarray}
\nonumber
\bar{\mathcal K}_4^{3,4}&=&-\frac{199 \left(p_1^2+q_1^2\right){}^2}{768 L^8}-\frac{911 \left(p_2^2+q_2^2\right) \left(p_1^2+q_1^2\right)}{384 \sqrt{10} L^8}-\frac{49 \left(p_2^2+q_2^2\right){}^2}{512 L^8}.
\end{eqnarray}
As the next step, we introduce action-angle variables $(I_1,I_2,\phi_1,\phi_2),$ which are defined by the transformation \eqref{accionangulo} and incorporate the terms associated to the action $L$ dropped in the previous steps and undo the time scalings, arriving at 
$$\tilde{\mathcal K}^{3,4}=-\frac{1}{2 L^2}+\frac{\eta ^8}{8 L^6}-\frac{\eta ^{9} (4 \sqrt{7} I_1-9 I_2)}{24 L^7}+\frac{\eta^{10}(1990 I_1^2+1822 \sqrt{10} I_1 I_2+735 I_2^2)}{1920 L^8}+O(\eta^{11}),$$
where $\eta=\varepsilon^{1/4}.$ Finally, we are in the position to apply Han, Li and Yi's Theorem~\ref{Hanliyi}, taking $a=3,$ $n=3,$ $m_1=8,$ $m_2=9,$ $m_3=10,$ $n_0=n_1=1,$ $n_2=n_3=3,$ $I^{n_0}=I^{n_1}=L,$ $I^{n_2}=I^{n_3}=(L,I_1,I_2),$ $\bar I^{n_0}=\bar I^{n_1}=L,$ $\bar I^{n_2}=\bar I^{n_3}=(I_1,I_2).$ The vector of frequencies is given by
$$
\Omega(I)=\left(\frac{1}{L^3},-\frac{3}{4 L^7},-\frac{\sqrt{7}}{6 L^7},\frac{3}{8 L^7},-\frac{1990 I_1+911 \sqrt{10} I_2}{960 L^8},-\frac{911 \sqrt{10} I_1+735 I_2}{960 L^8}\right).
$$
The corresponding matrix with dimension $6\times4$ whose columns are $\Omega(I),$ $\partial \Omega(I)/\partial L,$ $\partial\Omega(I)/\partial I_1$ and $\partial\Omega(I)/\partial I_2$ reads
\begin{equation}
M_{\Omega}(L,I_1,I_2)=\left[
\begin{array}{cccc}
 \frac{1}{L^3} & -\frac{3}{L^4} & 0 & 0 \\[1ex]
 -\frac{3}{4 L^7} & \frac{21}{4 L^8} & 0 & 0 \\
 -\frac{\sqrt{7}}{6 L^7} & \frac{7 \sqrt{7}}{6 L^8} & 0 & 0 \\[1ex]
 \frac{3}{8 L^7} & -\frac{21}{8 L^8} & 0 & 0 \\[1ex]
 -\frac{1990 I_1+911 \sqrt{10} I_2}{960 L^8} & \frac{1990 I_1+911 \sqrt{10} I_2}{120 L^9} & -\frac{199}{96 L^8} & -\frac{911}{96 \sqrt{10} L^8} \\[1ex]
 -\frac{911 \sqrt{10} I_1+735 I_2}{960 L^8} & \frac{911 \sqrt{10} I_1+735 I_2}{120 L^9} & -\frac{911}{96 \sqrt{10} L^8} & -\frac{49}{64 L^8} \\[1ex]
\end{array}
\right].
\end{equation}
Since the rank of $M_\Omega$ is equal to four, there exists a 3-tori surrounding each periodic solution generated by the isolates equilibria $\mathbf E_{3,4}.$ 
 \end{proof}

\section{Second reduced space}\label{sec:Reducido2}
In the Section~\ref{sec:Reducido1}, we applied a Lie-Deprit transformation to introduce an artificial symmetry up to second-order terms. This process removed the mean anomaly $\ell$. Additionally, the operation resulted in another symmetry, as the argument of periapis $g$ is absent in the normalized terms of the Hamiltonian \eqref{HamPromediado}. This added symmetry enables a further reduction related to the central action driven by $G$. The invariants corresponding to this action are expressed by the components of the angular momentum vector
$$
\beta_1=\frac{x_1+y_1}{2},\quad\beta_2=\frac{x_2+y_2}{2},\quad\beta_3=\frac{x_3+y_3}{2}.
$$
Let us consider $z=(\beta_1,\beta_2,\beta_3).$ The Poisson structure is given by $J( z)=\hat{z},$ more specifically,
\begin{equation}
J(\beta_1,\beta_2,\beta_3)=\left[
\begin{array}{ccc}
0 &-\beta_3 & \beta_2  \\
\beta_3 &0 & -\beta_1  \\
-\beta_2 &\beta_1 &0  \\
\end{array}
\right],
\end{equation} 
with the Casimir function $\mathcal C=\frac{1}{2}(\beta_1^2+\beta_2^2+\beta_3^2).$ The Hamiltonian function in terms of the twice-reduces space invariants is given by
\begin{equation}\label{eq:HamSecondReduced}
\mathcal H_{\ell,g}=-\dfrac{1}{2L^2}+\dfrac{\varepsilon^2}{2}\frac{-\beta_1^2+2 \sqrt{3} \beta_1 \beta_2+\beta_2^2+8 \beta_3^2}{24 L^3 \left(\beta_1^2+\beta_2^2+\beta_3^2\right)^{5/2}}+O(\varepsilon^3).
\end{equation}
Truncating the Hamiltonian $\mathcal H_{\ell,g} $ to second order and scaling the independent variable the equations of motion are:
\begin{equation}
\dot \beta_1=\beta_3 \left(\sqrt{3} \text{$\beta_1$}-7 \text{$\beta_2$}\right),\quad \dot \beta_2=\beta_3 \left(9 \text{$\beta_1$} -\sqrt{3} \text{$\beta_2$}\right),\quad \dot \beta_3=-\sqrt{3} \text{$\beta_1$}^2-2 \text{$\beta_1$} \text{$\beta_2$}+\sqrt{3} \text{$\beta_2$}^2.
\end{equation}
This vector field results in three isolated relative equilibria, which correspond to quasi-periodic orbits in the second-order truncated Hamiltonian.
\begin{equation}\label{eq:eq2red}
\mathbf E_{1,2}=(0,0,\pm G),\quad \mathbf E_{3,4}=\left(\pm\dfrac{\sqrt{3}G}{2},\mp\dfrac{G}{2},0\right),\quad \mathbf E_{5,6}=\left(\pm\dfrac{G}{2},\pm\dfrac{\sqrt{3}G}{2},0\right).
\end{equation}

\subsection{KAM Tori}
This section focuses on identifying the existence of three-dimensional KAM tori associated with the relative equilibria $\mathbf{E}_i$ for $i=1, \dots, 4$ in the twice-reduced space. To accomplish this, we reformulate the Hamiltonian \eqref{eq:HamSecondReduced} into a more suitable, manageable form. This reformulation involves a series of transformations aimed at simplifying the analysis and facilitating the detection of KAM tori.

For the relative equilibria $\mathbf E_{1,2},$ we introduce the following symplectic coordinates, centered at the equilibria
\begin{equation}\label{change}
a=\frac{\sqrt{2} \beta_1}{\sqrt{G\pm\beta_3}},\quad b=\pm\frac{\sqrt{2} \beta_2}{\sqrt{G\pm\beta_3}},
\end{equation}
with the inverse,
\begin{equation}\nonumber
\beta_1=\frac{1}{2} a \sqrt{4 G-b^2-a^2},\quad \beta_2=\pm\frac{1}{2} b \sqrt{4 G-b^2-a^2},\quad \beta_3=\pm\frac{1}{2} \left(2 G-b^2-a^2\right).
\end{equation}
The expression of the Hamiltonian \eqref{eq:HamSecondReduced} in both charts is given by
\begin{equation}\label{eq:casifin}
\mathcal K_{\ell,g}^{1,2}=-\dfrac{1}{L^2}+\dfrac{\varepsilon^2}{2}\frac{32 G^2-4 G \left(7 b^2\mp2 \sqrt{3} ab+9 a^2\right)+\left(b^2+a^2\right) \left(7 b^2\mp2 \sqrt{3} ab+9 a^2\right)}{96 G^5 L^3}+O(\varepsilon^3).
\end{equation}
Then, considering the following $\varepsilon^{-1/4}$ symplectic chart 
$$a=\varepsilon^{1/8}\frac{\sqrt{5}\bar a-10 \bar b}{\sqrt{2}\ 15^{3/4}},\quad b=\varepsilon^{1/8}\frac{3^{3/4} \bar a}{\sqrt{2} \sqrt[4]{5}},$$
on $\mathcal K_{\ell,g}^{1,2},$ we arrive at
$\bar{\mathcal K}_{\ell,g}^{1,2}=\varepsilon^{-1/4}\bar{\mathcal K}_1^{1,2}+\bar{\mathcal K}_2^{1,2}+\varepsilon^{1/4}\bar{\mathcal K}_4^{1,2}+O(\varepsilon ^{1/2}).
$
After applying the action-angle coordinates
$\bar a=\sqrt{2I}\cos(\phi),$ $\bar b=\sqrt{2I}\sin(\phi),$ average $\bar{\mathcal K}_4^{1,2}$ with respect to $\phi,$ incorporating the terms associated to the action $L$ dropped in the previous
steps and undoing the time scalings, we arrive at
\begin{equation}\nonumber
\tilde{\mathcal K}^{1,2}(I,\phi)=-\frac{1}{2 L^2}+\frac{\eta^8}{3 G^3 L^3}-\frac{\eta^9\sqrt{\frac{5}{3}} I}{2 G^4 L^3}+\frac{\eta^{10} I^2}{3 G^5 L^3}+O(\eta^{11}),
\end{equation}
where $\eta=\varepsilon^{1/4}.$ Applying  Han, Li and Yi's Theorem~\ref{Hanliyi}, and taking $a=3,$ $n=3,$ $m_1=8,$ $m_2=9,$ $m_3=10,$ $n_0=1,$ $n_1=2,$ $n_2=n_3=3,$ $I^{n_0}=L,$ $I^{n_1}=(L,G),$ $I^{n_2}=I^{n_3}=(L,G,I),$ $\bar I^{n_0}=L,$ $\bar I^{n_1}=G,$ $\bar I^{n_2}=\bar I^{n_3}=I.$  The vector of frequencies is given by
$$
\Omega(I)=\left( \frac{1}{L^3},-\frac{1}{G^4 L^3},-\frac{\sqrt{\frac{5}{3}}}{2 G^4 L^3},\frac{2 I}{3 G^5 L^3}\right).
$$
The corresponding column matrix $(\Omega(I),\partial \Omega(I)/\partial L,\partial \Omega(I)/\partial G,\partial \Omega(I)/\partial I)$ is

$$M_\Omega(L,G,I)=\left[
\begin{array}{cccc}
 \frac{1}{L^3} & -\frac{3}{L^4} & 0 & 0 \\[0.5ex]
 -\frac{1}{G^4 L^3} & \frac{3}{G^4 L^4} & \frac{4}{G^5 L^3} & 0 \\[0.5ex]
 -\frac{\sqrt{5/3}}{2 G^4 L^3} & \frac{\sqrt{15}}{2 G^4 L^4} & \frac{2 \sqrt{5/3}}{G^5 L^3} & 0 \\[0.5ex]
 \frac{2 I}{3 G^5 L^3} & -\frac{2 I}{G^5 L^4} & -\frac{10I}{3 G^6 L^3} & \frac{2}{3 G^5 L^3} \\[0.5ex]
\end{array}
\right].$$
Since the rank of $M_\Omega$ is equal to three, there exist a 3-tori surrounding each periodic solution
generated by the isolates equilibria $\mathbf E_{1,2}.$

For the equilibria $\mathbf E_{3,4},$ the treatment is similar. Consider the Hamiltonian function \eqref{eq:HamSecondReduced} and apply the local transformation \eqref{change} arriving at $\mathcal K_{\ell,g}^{1,2}$ given in \eqref{eq:casifin}. The equilibria $\mathbf E_{3,4}$ in the new variables become
$$ \mathbf E_{3,4}=\left(\pm\frac{\sqrt{3G}}{\sqrt{2}},-\frac{\sqrt{G}}{\sqrt{2}}\right).$$
We move the system so that $\mathbf{E}_{3,4}$ is positioned at the origin. The Hamiltonian function of the translated system becomes
\begin{eqnarray}
\nonumber \mathcal K_{\ell,g}^{3,4}&=&-\frac{1}{L^2} -\frac{\varepsilon^2}{2}\left(
\frac{19 p^2 q}{16 \sqrt{6} G^{9/2} L^3}-\frac{17 p^3}{48 \sqrt{2} G^{9/2} L^3}-\frac{25 p q^2}{48 \sqrt{2} G^{9/2} L^3}+\frac{19 q^3}{16 \sqrt{6} G^{9/2} L^3}\right.\\
\nonumber &&\qquad +\frac{7 p^4}{96 G^5 L^3}-\frac{p^3 q}{16 \sqrt{3} G^5 L^3}+\frac{p^2 q^2}{6 G^5 L^3}-\frac{p q^3}{16 \sqrt{3} G^5 L^3}+\frac{3 q^4}{32 G^5 L^3}\\
\nonumber &&\qquad\left.+\frac{13 p^2}{48 G^4 L^3}-\frac{3 \sqrt{3} p q}{8 G^4 L^3}+\frac{31 q^2}{48 G^4 L^3}-\frac{1}{12 G^3 L^3}
\right)+O(\varepsilon^3).
\end{eqnarray}
Using the following $\varepsilon^{-1/4}$ symplectic chart 
$$a= \varepsilon^{1/8}\frac{9 \sqrt{30} \bar a-40 \bar b}{2\ 10^{3/4} \sqrt{31}},\quad b=\varepsilon^{1/8}\frac{\sqrt{31} \bar a}{2 \sqrt[4]{10}},$$
on $\mathcal K_{\ell,g}^{3,4},$ we obtain
$\bar{\mathcal K}_{\ell,g}^{3,4}=\varepsilon^{-1/4}\bar{\mathcal K}_1^{3,4}+\bar{\mathcal K}_2^{3,4}+\varepsilon^{1/8} \bar{\mathcal K}_3^{3,4}+\varepsilon^{1/4}\bar{\mathcal K}_4^{3,4}+O(\varepsilon ^{1/2}).
$
After applying the action-angle coordinates
$\bar a=\sqrt{2I}\cos(\phi),$ $\bar b=\sqrt{2I}\sin(\phi),$ averaging $\bar{\mathcal K}_3^{3,4}$ and $\bar{\mathcal K}_4^{3,4}$ with respect to $\phi,$ incorporating the terms associates to the action $L$ dropped in the previous
steps and undoing the time scalings, we arrive at
\begin{equation}\nonumber
\tilde{\mathcal K}^{3,4}(I,\phi)=-\frac{1}{2 L^2}-\frac{\eta^8}{12 G^3 L^3}+\frac{\eta^9\sqrt{5/2} I}{3 G^4 L^3}+\frac{\eta^{10}199 I^2}{192 G^5 L^3}+O(\eta^{11})
\end{equation}
where $\eta=\varepsilon^{1/4}.$ Considering $a=3,$ $n=3,$ $m_1=8,$ $m_2=9,$ $m_3=10,$ $n_0=1,$ $n_1=2,$ $n_2=n_3=3,$ $I^{n_0}=L,$ $I^{n_1}=(L,G),$ $I^{n_2}=I^{n_3}=(L,G,I),$ $\bar I^{n_0}=L,$ $\bar I^{n_1}=G,$ $\bar I^{n_2}=\bar I^{n_3}=I,$ we apply  Han, Li and Yi's Theorem~\ref{Hanliyi}. The vector of frequencies is given by
$$
\Omega(I)=\left( \frac{1}{L^3},\frac{1}{4 G^4 L^3},\frac{\sqrt{5/2}}{3 G^4 L^3},\frac{199 I}{96 G^5 L^3}\right).
$$
The corresponding column matrix $(\Omega(I),\partial \Omega(I)/\partial L,\partial \Omega(I)/\partial G,\partial \Omega(I)/\partial I)$ is
$$M_\Omega=\left[
\begin{array}{cccc}
 \frac{1}{L^3} & -\frac{3}{L^4} & 0 & 0 \\[1ex]
 \frac{1}{4 G^4 L^3} & -\frac{3}{4 G^4 L^4} & -\frac{1}{G^5 L^3} & 0 \\[1ex]
 \frac{\sqrt{5/2}}{3 G^4 L^3} & -\frac{\sqrt{5/2}}{G^4 L^4} & -\frac{2 \sqrt{10}}{3 G^5 L^3} & 0 \\[1ex]
 \frac{199 I}{96 G^5 L^3} & -\frac{199 I}{32 G^5 L^4} & -\frac{995 I}{96 G^6 L^3} & \frac{199}{96 G^5 L^3} \\[1ex]
\end{array}
\right].$$
Since the rank of $M_\Omega$ is equal to three, there exists a 3-tori surrounding each periodic solution
generated by the isolates equilibria $\mathbf E_{3,4}.$ We have proven the following result.

\begin{Theorem}
There are families of invariant KAM 3-tori of the full system around the cuasi-periodic solutions $\mathbf E_{1,2}$ and $\mathbf E_{3,4}$ given in \eqref{eq:eq2red}. These invariant tori form a majority in the sense that the measure of the complement of their union is of the order $O(\varepsilon^{\delta/4}),$ with $0<\delta<\frac{1}{5}.$
\end{Theorem}

\appendix

\section{Reeb's Theorem and Parametrical Stable}\label{apA}
In this appendix we recall to the reader some result associated with Reeb's Theorem taken from \cite{Palacian2008}. Let us consider the linear constant coefficient Hamiltonian system 
\begin{equation}\label{eq:coro}
\dot y=Cy=J\nabla H(y),\qquad H=\dfrac{1}{2}y^TSy,
\end{equation}
where $S$ is a symmetric matrix and $C=JS$ is a Hamiltonian matrix. 
\begin{Definition}[parametrically stability]
The system \eqref{eq:coro} is parametrically stable if it and all sufficiently small linear constant Hamiltonian perturbations of it are stable. If the system \eqref{eq:coro} is stable but not strongly stable, we shall say that it is weakly stable.
\end{Definition}
Let $\pm \alpha_1i,\pm\alpha_2 i,\dots,\pm\alpha_si$ be the eigenvalues of the stable matrix $C,$ and $V_j,$ $j=1,\dots,s,$ be the maximal real linear subspace where $C$ has eigenvalues $\pm\alpha_ji$. So $V_j$ is a $C$-invariant symplectic subspace, $C$ restricted to $V_j$ has eigenvalues $\pm\alpha_ji,$ and $\mathbb R^{2n}=V_1\oplus V_2\oplus\dots\oplus V_s.$ Let $H_j$ be the restriction of $H$ to $V_j.$
\begin{Theorem}[Krein-Gel'fand]
System \eqref{eq:coro} is parametrically stable if and only if 
\begin{itemize}
\item all the eigenvalues of $C$ are purely imaginary,
\item $C$ is nonsingular,
\item $C$ is diagonalizable over the complex numbers, and
\item the Hamiltonian $H_j$ is positive or negative definite for each $j.$
\end{itemize}
\end{Theorem}
This proof can also be found in \cite{me,yaku}.
Let $\varepsilon$ be a small parameter, $\mathcal H_1:M\to\mathbb R$ be smooth, $\mathcal H_\varepsilon=\mathcal H_0+\varepsilon \mathcal H_1,$ $Y_\varepsilon=Y_0+\varepsilon Y_1=d\mathcal H_\varepsilon^\#,$ $\mathcal N_\varepsilon(h)=\mathcal H_\varepsilon^{-1}(h),$ $\pi:\mathcal N_\varepsilon(h)\to\mathcal B(h)$ be the proyection and $\phi_\varepsilon^t$ be the flow defined by $Y_\varepsilon.$ 
Let the average of $\mathcal H_1$ be
\begin{equation}
\bar{\mathcal H}=\frac{1}{T} \int_0^T\mathcal H_1(\phi_0^t)dt.
\end{equation}

\begin{Theorem}[Reeb]\label{teo:Reeb} If $\bar{\mathcal H}$ has a nondegenerate critical point at $\pi(p)=\bar p\in \mathcal B(h)$ with $p\in \mathcal N_0,$ then there are smooth functions $p(\varepsilon)$ and $T(\varepsilon)$ for $\varepsilon$ small with $p(0)=p,$ $T(0)=T,$ and $p(\epsilon)\in\mathcal N_\varepsilon$ and the solution of $Y_\varepsilon$ through $p(\epsilon)$ is $T(\varepsilon)$-periodic. In addition, if the characteristic exponents of the critical point $\bar p$ (that is, the eigenvalues of the matrix $A=JD^2 \bar{ \mathcal H }(\bar p))$ are $\lambda_1,\lambda_2,\dots,\lambda_{2n-2},$ then the characteristic multipliers of the periodic solution through $p(\varepsilon)$ are
$$1,1,1+\varepsilon \lambda_1 T+O(\varepsilon^2), 1+\varepsilon \lambda_2 T+O(\varepsilon^2),\dots,1+\varepsilon \lambda_{2n-2} T+O(\varepsilon^2).$$
\end{Theorem}

\begin{Theorem}\label{Teo:inestable}
Let $p$ and $\bar p$ be as in the previous Theorem. If one or more of the characteristic exponents $\lambda_j$ is real or has nonzero real part, then the periodic solution through $p(\varepsilon)$ is unestable. If the matrix $A$ is parametrically stable, then the periodic solution through $p(\varepsilon)$ is elliptic, i.e., linearly stable.
\end{Theorem}

For more information on this subject the reader is referred to \cite{Reeb1952,Palacian2008,palacian2011}.

\section{KAM Theory: Han Li and Yi's Theorem}\label{apB}

This section relies in \cite{hanliyi} and collect the results from KAM theory that we employ to treat a Hamiltonian with high degeneracy. We follow the notation given in \cite{meyerLunar} starting with a Hamiltonian system of the form 
\begin{equation}\label{eq:24}
H(I,\varphi,\varepsilon)=h_0(I^{n_0})+\varepsilon^{m_1}h_1(I^{n_1})+\dots+\varepsilon^{m_a}h_a(I^{n_a})+\varepsilon^{m_a+1}P(I,\varphi,\varepsilon),
\end{equation}
where $(I,\varphi)\in\mathbb R^n\times\mathbb T^n$ are the action-angle variables with the standard symplectic structure $dI\wedge d\varphi,$ and $\varepsilon >0$ is a sufficiently small parameter. The Hamiltonian $H$ is real analytic in
$(I,\varphi,\varepsilon)$ and in particular $p$ is a smooth in $\varepsilon$. The parameters $a,m,m_i$ $(i= 0,1,\dots,a)$ and
$m_j$ $(j=1,2,\dots,a)$, are positive integers satisfying $n_0 \leq n_1 \leq \dots\leq n_a = n,$ $m_1 \leq m_2 \leq \dots \leq
m_a = m,$ $I^{n_i} = (I_1,\dots ,I^{n_i} )$, for $i= 1,2,\dots,a.$

The Hamiltonian $H(I,\varphi,\varepsilon)$ is considered in a bounded closed region $Z\times \mathbb T^n  \times [0,\varepsilon^*] \subset
\mathbb R^n \times \mathbb T^n \times [0,\varepsilon^*]$ for some fixed $\varepsilon^*$ with $0 <\varepsilon^*<1.$ For each $\varepsilon$ the integrable part of $H,$
$$X_\varepsilon(I)=h_0(I^{n_0})+\varepsilon^{m_1}h_1(I^{n_1})+\dots +\varepsilon^{m_1}h_a(I^{n_a}),$$
admits a family of invariant $n$-tori $T^\varepsilon_\zeta= \{\zeta\}\times\mathbb T^n$ with linear flows $x_0 + \omega^\varepsilon (\zeta)t,$ where, for
each $\zeta \in Z,$ $\omega^{\varepsilon}(\zeta) = \nabla X_\varepsilon (\zeta)$ is the frequency vector of the $n$-torus $T^\varepsilon_\zeta$ and $\nabla$ is the gradient operator. When $\omega^\varepsilon(\zeta)$ is non-resonant, the flow on the $n$-torus $T^\varepsilon_\zeta$ becomes quasi-periodic with slow and fast frequencies of different scales. We refer to the integrable part $X_\varepsilon$ and its associated tori $\{T^\varepsilon_\zeta\}$ as the intermediate Hamiltonian and intermediate tori, respectively.

Let $\bar I^{n_i}=(I_{n_{i-1}+1},\dots,I_{n_i}),$ $i=0,1,\dots,a$ (where $n_{i-1}=0,$ hence $\bar I^{n_0}=I^{n_0}$), and define $$\Omega=(\nabla_{\bar I^{n_0}}h_0(I^{n_0}),\dots,\nabla_{\bar I^{n_a}}h_{n_a}(I^{n_a})),$$
such that for each $i=0,1,\dots,a,$ $\nabla_{\bar I^{n_i}}$ denotes the gradient with respect to $\bar I^{n_i}.$

The following theorem gives the right setting in which one can ensure the persistence of KAM tori for the Hamiltonian \eqref{eq:24}

\begin{Theorem}[Han, Li and Yi]\label{Hanliyi}
Let $\delta$ be given with $0<\delta<1/5.$ Assume that there is a positive integer s such that
$$Rank\{\partial_I^\alpha \Omega(I):0\leq \vert \alpha\vert\leq s\}=n,\quad \forall I\in Z.$$
Then, there exists an $\varepsilon_0>0$ and a family of Cantor sets $Z_\varepsilon\subset Z,$ $0<\varepsilon<\varepsilon_0,$ such that each $\zeta$ corresponds to a real analitic, invariant, quasi-periodic n-torus $T_\zeta^\varepsilon.$ The measure of $Z\setminus Z_\varepsilon$ is $O(\varepsilon ^{\delta/s})$ and the family $\{\bar T_\zeta^{\varepsilon}:\zeta \in Z_\varepsilon,0<\varepsilon<\varepsilon_0\}$ varies Whitney smoothly.
\end{Theorem}

\section{Complete Expression of $\mathcal H_2^0$}
{\footnotesize \begin{eqnarray}
\nonumber \mathcal H_2^0&=&\frac{\left(c_E-e\right){}^2 c_g^2 H^2}{2 G^2 L^6 \left(1-e c_E\right){}^5}+\frac{c_h^2 H^2}{8 G^2 L^6 \left(1-e c_E\right){}^3}+\frac{\left(1-e^2\right) c_g^2 c_h^2 s_E^2 H^2}{8 G^2 L^6 \left(1-e c_E\right){}^5}+\frac{\left(c_E-e\right){}^2 c_h^2 s_g^2 H^2}{8 G^2 L^6 \left(1-e c_E\right){}^5}\\
\nonumber &&
+\frac{\left(1-e^2\right) s_E^2 s_g^2 H^2}{2 G^2 L^6 \left(1-e c_E\right){}^5}+\frac{\left(c_E-e\right){}^2 c_g^2 s_h^2 H^2}{8 G^2 L^6 \left(1-e c_E\right){}^5}+\frac{\left(1-e^2\right) s_E^2 s_g^2 s_h^2 H^2}{8 G^2 L^6 \left(1-e c_E\right){}^5}+\frac{\eta  \left(c_E-e\right) c_g c_h^2 s_E s_g H^2}{2 G^2 L^6 \left(1-e c_E\right){}^5}\\
\nonumber&&
+\frac{\left(1-e^2\right) \sqrt{3} c_g^2 c_h s_E^2 s_h H^2}{4 G^2 L^6 \left(1-e c_E\right){}^5}+\frac{\sqrt{3} \left(c_E-e\right){}^2 c_h s_g^2 s_h H^2}{4 G^2 L^6 \left(1-e c_E\right){}^5}+\frac{\sqrt{3} c_h s_h H^2-2H^2}{4 G^2 L^6 \left(1-e c_E\right){}^3}-\frac{s_h^2 H^2}{8 G^2 L^6 \left(1-e c_E\right){}^3}\\
\nonumber&&
+\frac{\eta  \sqrt{3} \left(c_E-e\right) c_g c_h s_E s_g s_h H^2}{G^2 L^6 \left(1-e c_E\right){}^5}-\frac{2 \eta  \left(c_E-e\right) c_g s_E s_g H^2}{G^2 L^6 \left(1-e c_E\right){}^5}-\frac{\left(1-e^2\right) c_g^2 s_E^2 H^2}{2 G^2 L^6 \left(1-e c_E\right){}^5}-\frac{\left(c_E-e\right){}^2 s_g^2 H^2}{2 G^2 L^6 \left(1-e c_E\right){}^5}\\
\nonumber &&
-\frac{\eta  \left(c_E-e\right) c_g s_E s_g s_h^2 H^2}{2 G^2 L^6 \left(1-e c_E\right){}^5}-\frac{\sqrt{3} \left(1-e^2\right) c_h s_E^2 s_g^2 s_h H^2}{4 G^2 L^6 \left(1-e c_E\right){}^5}-\frac{\sqrt{3} \left(c_E-e\right){}^2 c_g^2 c_h s_h H^2}{4 G^2 L^6 \left(1-e c_E\right){}^5}\\
\nonumber &&
-\frac{\left(c_E-e\right){}^2 c_g^2 c_h^2 H^2}{8 G^2 L^6 \left(1-e c_E\right){}^5}-\frac{\left(1-e^2\right) c_h^2 s_E^2 s_g^2 H^2}{8 G^2 L^6 \left(1-e c_E\right){}^5}-\frac{\left(1-e^2\right) c_g^2 s_E^2 s_h^2 H^2}{8 G^2 L^6 \left(1-e c_E\right){}^5}-\frac{\left(c_E-e\right){}^2 s_g^2 s_h^2 H^2}{8 G^2 L^6 \left(1-e c_E\right){}^5}\\
\nonumber &&
+\frac{\eta  \sqrt{3} \left(c_E-e\right) c_h^2 s_E s_g^2 H}{2 G L^6 \left(1-e c_E\right){}^5}+\frac{\eta  \sqrt{3} \left(c_E-e\right) c_g^2 s_E s_h^2 H}{2 G L^6 \left(1-e c_E\right){}^5}+\frac{\sqrt{3} \left(c_E-e\right){}^2 c_g s_g s_h^2 H}{2 G L^6 \left(1-e c_E\right){}^5}\\
\nonumber&&
+\frac{\left(1-e^2\right) \sqrt{3} c_g c_h^2 s_E^2 s_g H}{2 G L^6 \left(1-e c_E\right){}^5}+\frac{\eta  \left(c_E-e\right) c_g^2 c_h s_E s_h H}{G L^6 \left(1-e c_E\right){}^5}+\frac{\left(c_E-e\right){}^2 c_g c_h s_g s_h H}{G L^6 \left(1-e c_E\right){}^5}\\
\nonumber &&
-\frac{\eta  \left(c_E-e\right) c_h s_E s_g^2 s_h H}{G L^6 \left(1-e c_E\right){}^5}-\frac{\left(1-e^2\right) c_g c_h s_E^2 s_g s_h H}{G L^6 \left(1-e c_E\right){}^5}-\frac{\sqrt{3} \eta  \left(c_E-e\right) s_E s_g^2 s_h^2 H}{2 G L^6 \left(1-e c_E\right){}^5}-\\
\nonumber &&
\frac{\sqrt{3} \left(1-e^2\right) c_g s_E^2 s_g s_h^2 H}{2 G L^6 \left(1-e c_E\right){}^5}-\frac{\sqrt{3} \eta  \left(c_E-e\right) c_g^2 c_h^2 s_E H}{2 G L^6 \left(1-e c_E\right){}^5}-\frac{\sqrt{3} \left(c_E-e\right){}^2 c_g c_h^2 s_g H}{2 G L^6 \left(1-e c_E\right){}^5}\\
\nonumber &&
+\frac{\left(1-e^2\right) c_g^2 s_E^2}{2 L^6 \left(1-e c_E\right){}^5}+\frac{\left(1-e^2\right) c_g^2 c_h^2 s_E^2}{8 L^6 \left(1-e c_E\right){}^5}+\frac{\left(c_E-e\right){}^2 c_h^2 s_g^2}{8 L^6 \left(1-e c_E\right){}^5}+\frac{\left(c_E-e\right){}^2 s_g^2}{2 L^6 \left(1-e c_E\right){}^5}\\
\nonumber&&
+\frac{s_h^2}{8 L^6 \left(1-e c_E\right){}^3}+\frac{\eta  \left(c_E-e\right) c_g c_h^2 s_E s_g}{2 L^6 \left(1-e c_E\right){}^5}+\frac{2 \eta  \left(c_E-e\right) c_g s_E s_g}{L^6 \left(1-e c_E\right){}^5}+\frac{\left(1-e^2\right) \sqrt{3} c_g^2 c_h s_E^2 s_h}{4 L^6 \left(1-e c_E\right){}^5}\\
\nonumber &&
+\frac{\sqrt{3} \left(c_E-e\right){}^2 c_h s_g^2 s_h}{4 L^6 \left(1-e c_E\right){}^5}+\frac{\eta  \sqrt{3} \left(c_E-e\right) c_g c_h s_E s_g s_h}{L^6 \left(1-e c_E\right){}^5}-\frac{\sqrt{3} c_h s_h}{4 L^6 \left(1-e c_E\right){}^3}+\frac{1}{6 L^6 \left(1-e c_E\right){}^3}\\
\nonumber &&-\frac{c_h^2}{8 L^6 \left(1-e c_E\right){}^3}-\frac{\left(c_E-e\right){}^2 c_g^2}{2 L^6 \left(1-e c_E\right){}^5}-\frac{\left(1-e^2\right) s_E^2 s_g^2}{2 L^6 \left(1-e c_E\right){}^5}-\frac{\eta  \left(c_E-e\right) c_g s_E s_g s_h^2}{2 L^6 \left(1-e c_E\right){}^5}\\
\nonumber &&-\frac{\sqrt{3} \left(1-e^2\right) c_h s_E^2 s_g^2 s_h}{4 L^6 \left(1-e c_E\right){}^5}-\frac{\sqrt{3} \left(c_E-e\right){}^2 c_g^2 c_h s_h}{4 L^6 \left(1-e c_E\right){}^5}-\frac{\left(c_E-e\right){}^2 c_g^2 c_h^2}{8 L^6 \left(1-e c_E\right){}^5}-\frac{\left(1-e^2\right) c_h^2 s_E^2 s_g^2}{8 L^6 \left(1-e c_E\right){}^5}\\
\nonumber &&
-\frac{\left(1-e^2\right) c_g^2 s_E^2 s_h^2}{8 L^6 \left(1-e c_E\right){}^5}-\frac{\left(c_E-e\right){}^2 s_g^2 s_h^2}{8 L^6 \left(1-e c_E\right){}^5}.
\end{eqnarray}
}

\section*{Declarations}
\subsection*{Ethical approval} Not applicable.
\subsection*{Consent to participate} Not applicable.
\subsection*{Consent to publish} All authors have provided their consent for the publication of this article.
\subsection*{Funding}The authors ET and FC has been partially funded by Research Agencies of Chile and Universidad del Bío-Bío is acknowledged, they came in the form of research projects GI2310532-VRIP-UBB and ANID PhD/2021-21210522. The third author was partially supported by ANID-Chile through FONDECYT Iniciación 11240582.

\bibliographystyle{abbrv}
\bibliography{Bibliography}
\end{document}